\documentclass[a4paper,12pt]{amsproc}
\usepackage{amsmath,amsfonts,amssymb,amsthm,anysize}
\usepackage{enumerate}
\usepackage{epsfig}
\usepackage{supertabular}
\usepackage{graphicx}

\newcommand{\tr}{\text{Tr}}

\newcommand\Z{{\mathbb Z}}

\newcommand\F{{\mathbb F}}

\newcommand\Tr{{\mathrm{Tr}}}
\newcommand\Norm{{\mathrm{Norm}}}
\newcommand\Cay{{\mathrm{Cay}}}

\newcommand\PG{\mathsf{PG}}
\renewcommand\mod{{\mathrm{mod\, \, }}}

\newcommand\cQ{{\mathcal Q}}

\newcommand{\dualG}{{\widehat G}}

\theoremstyle{plain}
\newtheorem{theorem}{Theorem}[section]
\newtheorem{problem}[theorem]{Problem}
\newtheorem{lemma}[theorem]{Lemma}
\newtheorem{corollary}[theorem]{Corollary}

\newtheorem{proposition}[theorem]{Proposition}
\newtheorem{example}[theorem]{Example}

\newtheorem{result}[theorem]{Result}
\numberwithin{equation}{section}

\theoremstyle{remark}
\newtheorem{remark}[theorem]{Remark}

\usepackage{url, hyperref}
\hypersetup{citecolor=blue, linkcolor=blue, colorlinks=true}

\renewcommand\le{\leqslant}

\begin{document}

\title[Strongly regular Cayley graphs]{Construction of strongly regular Cayley graphs based on three-valued Gauss periods}


\author{Koji Momihara}
\address{ %
Department of Mathematics\\
Faculty of Education\\
Kumamoto University\\
2-40-1 Kurokami, Kumamoto 860-8555, Japan}
\email{momihara@educ.kumamoto-u.ac.jp}
\thanks{The author acknowledges the support by 
JSPS under Grant-in-Aid for Young Scientists (B) 17K14236 and Scientific Research (B) 15H03636.}


\subjclass[2010]{05E30, 11T22 (primary), 05C50, 05B10, 11T24 (secondary)}
\keywords{Gauss period, Gauss sum, Cayley graph, strongly regular graph, cyclotomic class}

\begin{abstract}
In this paper, we give a construction of strongly regular Cayley graphs on the additive groups of finite fields 
based on three-valued Gauss periods. 
As consequences, we obtain two infinite families and one  sporadic example of new strongly regular Cayley graphs. This construction  can be viewed as a generalization of that of strongly regular Cayley graphs obtained in \cite{BLMX}.  
\end{abstract}


\maketitle

\section{Introduction}\label{sec:intro}

A {\it strongly regular graph} with parameters $(v,k,\lambda,\mu)$ is a simple undirected $k$-regular graph $\Gamma$ on the set of $v$ vertices satisfying the following conditions: 
for any two vertices $x$ and $y$, 
\[
|\{z\in V(\Gamma):(x,z),(y,z)\in E(\Gamma)\}|=
\left\{
\begin{array}{ll}
\lambda,&  \mbox{if $(x,y)\in E(\Gamma)$,}\\
\mu,&  \mbox{if $(x,y)\not\in E(\Gamma)$,}
 \end{array}
\right.
\]
where $V(\Gamma)$ is the set of vertices of $\Gamma$ and 
$E(\Gamma)$ is the set of edges of $\Gamma$. 
A strongly regular graph is said to be of {\it Latin square type} (resp. {\it  negative Latin square type}) if it has parameters 
$(v,k,\lambda,\mu)=(n^2,r(n-\epsilon),\epsilon n+r^2-3\epsilon r,r^2-\epsilon r)$ with $\epsilon=1$ (resp. $\epsilon=-1$).  
A regular graph, not complete or edgeless, is strongly regular if and only if its adjacency matrix has exactly two restricted eigenvalues~\cite[Theorem 9.1.2]{BH}. Here,  we say that an eigenvalue of the adjacency matrix is {\it restricted} if it has an eigenvector perpendicular to the all-ones vector. 

An effective way for constructing strongly regular graphs is to use Cayley graphs.  
Let $G$ be a finite (additively written) abelian group and $D$ be an inverse-closed subset of $G\setminus \{0_G\}$, where $0_G$ is the identity of $G$.  We define a graph $\Cay(G,D)$ with the elements of $G$ as its vertices; two vertices $x$ and $y$ are adjacent if and only if $x-y\in D$. The graph $\Cay(G,D)$ is called a {\it Cayley graph} on $G$ with connection set $D$. 
The eigenvalues of $\Cay(G,D)$ are given by $\psi(D)$, $\psi\in \dualG$, where $\dualG$ is the group consisting of all characters of $G$. By the aforementioned characterization of strongly regular graphs, $\Cay(G,D)$ with connection set $D$($\not=\emptyset,G$) is strongly regular if and only if $\psi(D)$, $\psi\in \dualG\setminus\{\psi_0\}$, take exactly two values, where $\psi_0$ is the identity of $\dualG$.  
If $\Cay(G,D)$ is strongly regular, the  set $D$ is often called a {\it partial difference set} on $G$. For basic results on strongly regular Cayley graphs and partial difference sets, we refer the reader to the survey \cite{Ma}. 

A classical method for constructing strongly regular Cayley graphs in the 
additive groups (i.e., elementary abelian groups) of finite fields is to use 
cyclotomic classes. 
Let $q$ be a prime power and $\F_q$ be the finite field of order $q$. Furthermore, 
let $\omega$ be a fixed primitive element of $\F_q$ and  $N$ be a positive integer dividing $q-1$. For $0\le i\le N-1$ we set $C_i^{(N,q)}=\omega^i C_0$, where $C_0$ is the subgroup of order $(q-1)/N$ of $\F_q^\ast$. Here, $\F_q^\ast$ denotes the multiplicative group of $\F_q$. These cosets are called {\it cyclotomic classes} 
of order $N$ of $\F_q$. One may take a union of some cyclotomic classes of order $N$ of $\F_q$ as the connection set of a Cayley graph. Many researchers have studies the problem of determining when a Cayley graph with a union of cyclotomic classes 
as its connection set forms a strongly regular graph. However, this method has had only limited success. One of the 
reason why not so many strongly regular graphs have been discovered 
by this method is the difficulty of evaluating Gauss sums with respect to multiplicative characters 
of  large order.  In fact, by the orthogonality of characters, the character values of a union of cyclotomic classes of order $N$ can be represented as 
sums of linear combinations of Gauss sums with respect to multiplicative characters 
of  exponent $N$. 

Typical cases where the Gauss sums of order $N$ of $\F_q$ have been evaluated are listed below: 
\begin{itemize}
\item[(1)] (small order case~\cite{BEW97}) $N\le 24$ (but some of the evaluations are not explicit),  
\item[(2)] (semi-primitive case~\cite{BEW97}) $-1\in \langle p\rangle \,(\mod{N})$, 
\item[(3)] (index $2$ case~\cite{YX}) $[\Z_N^\times:\langle p\rangle]=2$, where $\Z_N^\times$ is the unit group of $\Z_N$, 
\end{itemize}
where $p$ is the characteristic of $\F_{q}$. 
Corresponding to these cases, some constructions of strongly regular Cayley graphs have been given. For example, if $N=2$ and $q\equiv 1\,(\mod{4})$, $\Cay(\F_q,C_0^{(N,q)})$ forms a strongly regular graph, the so-called {\it Paley graph}. Brouwer, Wilson and Xiang~\cite{BWX}  constructed strongly regular 
Cayley graphs based on cyclotomic classes in semi-primitive case.  
Recently, Feng and Xiang~\cite{FX}  gave a construction of strongly regular Cayley graphs based on cyclotomic classes in index $2$ case. 
For related results, we refer the reader to \cite{CK,FMX3,GXY,LS,M} and references therein. 

On the other hand, strongly regular Cayley graphs on finite fields have been 
studied in finite geometry. In particular, finite geometric objects, such as {\it $m$-ovoids} and {\it $i$-tight sets}, can give arise to strongly regular Cayley graphs on the additive group of finite fields~\cite{BKLP,BLP,BD,CP,D,FMX1,L1,L2,S}. In particular, 
see \cite[Section~6]{BKLP} for the relationship between these finite geometric objects and strongly regular graphs. In these studies, instead of evaluating Gauss sums,  the character values were directly computed by geometric or  
group-theoretic 
arguments. Recently, Bamberg, Lee, Xiang and the author~\cite{BLMX}  constructed strongly regular Cayley graphs corresponding to  $\frac{q+1}{2}$-ovoids in the elliptic quadric ${\mathcal Q}^-(5,q)$, which are  also corresponding to a finite geometric object, so-called {\it hemisystems} of the Hermitian surface. In particular, the authors used cyclotomic classes of order $4(q^2+q+1)$ of $\F_{q^6}$ combining with a geometric argument. As a more theoretical understanding,  Xiang and the author~\cite{MX} realized that two-valued Gauss periods and a partition of the Singer difference set are behind the construction. Furthermore, the authors showed that the construction can be done in a more general setting  within the framework of two-valued Gauss periods.   

In this paper, we show that the construction given in \cite{BLMX} can also work within the framework of three-valued Gauss periods. As a consequence, we obtain the following two infinite families of new strongly regular Cayley graphs. 
\begin{theorem}
\label{thm:result}
There exists a strongly regular Cayley graph on $(\F_{q^6},+)$ with negative Latin square type parameters $(q^6,r(q^3+1),q^3+r^2-3r,r^2-r)$, where  
$r=M(q^2-1)/2$, in the following cases: 
\begin{itemize}
\item[(i)] $M=3$ and $q\equiv 7\,(\mod{24})$, 
\item[(ii)] $M=7$ and  $q\equiv 11,51\,(\mod{56})$. 
\end{itemize}
\end{theorem}
In the case where $M=1$ and $q\equiv 3\,(\mod{4})$ in the theorem above, we can also obtain strongly regular Cayley graphs, which were already found in \cite{BLMX}.  
Thus, our construction can be viewed as a generalization of that given in \cite{BLMX}. 

This paper is organized as follows. In Section~\ref{section:prelim}, we 
give the background
on Gauss sums, Gauss periods and Cayley graphs on finite fields. In Section~\ref{section:infinite}, we review known results on three-valued Gauss periods and find two families of three-valued Gauss periods forming an arithmetic progression.  In Section~\ref{section:construction}, we give a construction of Cayley graphs based on three-valued Gauss periods. Furthermore, we give a sufficient condition for the Cayley graph to be strongly regular. 
In Section~\ref{sec:find}, we  show that some of the three-valued Gauss periods found in Section~\ref{section:infinite} satisfy the sufficient condition. Then, we obtain our main result. 
In the final section,  we give one 
 sporadic example of a strongly regular Cayley graph on $(\F_{7^{14}},+)$ not belonging to the two families of strongly regular graphs in Theorem~\ref{thm:result}. 
\section{Preliminaries}\label{section:prelim}
\subsection{Gauss sums}
In this section, we collect basic results on Gauss sums.  

For a multiplicative character
$\chi$  and the canonical
additive character $\psi_{\F_q}$ of $\F_q$, define the {\it Gauss sum} by
\[
G_q(\chi)=\sum_{x\in \F_q^\ast}\chi(x)\psi_{\F_q}(x).
\]
If $\chi$ is of order $N$, we may say that the Gauss sum is of order $N$. 
We will use the following facts without preamble. 
\begin{enumerate}
\item[(i)] $G_q(\chi)\overline{G_q(\chi)}=q$ if $\chi$ is nontrivial;
\item[(ii)] $G_q(\chi^{-1})=\chi(-1)\overline{G_q(\chi)}$;
\item[(iii)] $G_q(\chi)=-1$ if $\chi$ is trivial.
\end{enumerate}
The Gauss sum with respect to the quadratic character is explicitly computable as follows. 
\begin{theorem}\cite[Theorem~5.15]{LN97} \label{thm:Gauss}
Let $q=p^f$ be a prime power with $p$ a prime and $\eta$ be the quadratic character of $\F_q$. 
Then, 
\begin{equation}\label{eq:Gaussquad1}
G_q(\eta)=\begin{cases}
(-1)^{f-1}q^{1/2}& \text{ if }  p\equiv 1\,(\mod{4}), \\
(-1)^{f-1}i^f q^{1/2} & \text{ if } p\equiv 3\,(\mod{4}). 
\end{cases}
\end{equation}
\end{theorem}
Also, in semi-primitive case, the 
Gauss sum is computable. 
\begin{theorem}\label{thm:semiprim}{\em (\cite[Theorem~11.6.3]{BEW97})}
Let $p$ be a prime. 
Suppose that $N>2$ and $p$ is semi-primitive modulo $N$, 
i.e., there exists a positive integer $j$ such that  $p^j\equiv -1\pmod{N}$. Choose 
$j$ minimal and write 
$f=2js$ for any positive integer $r$. Let $\chi$ be a multiplicative character of order $N$ of $\F_{p^f}$. 
Then, 
\[
p^{-f/2}G_{p^f}(\chi)=
\left\{
\begin{array}{ll}
(-1)^{s-1},&  \mbox{if $p=2$,}\\
(-1)^{s-1+(p^j+1)s/N},&  \mbox{if $p>2$. }
 \end{array}
\right.
\]
\end{theorem}
The following is often referred to as the {\it Davenport-Hasse lifting formula}. 
\begin{theorem}\label{thm:lift}
{\em (\cite[Theorem~11.5.2]{BEW97})}
Let $\chi'$ be a nontrivial multiplicative character of $\F_{q}$ and 
let $\chi$ be the lift of $\chi'$ to $\F_{q^{m}}$, i.e., $\chi(\alpha)=\chi'(\Norm_{q^{m}/q}(\alpha))$ for $\alpha\in \F_{q^{m}}$, where $m\geq 2$ is an integer. Then 
\[
G_{q^m}(\chi)=(-1)^{m-1}(G_{q}(\chi'))^m. 
\]
\end{theorem}
Also, we need the following formula on Gauss sums. 
\begin{theorem}\label{Gaussmain}{\em (\cite[Corollary~2.8]{BLMX})}
Let 
$q^m\equiv 3\,(\mod{4})$ be an odd prime power and $N$ be an odd positive integer dividing $(q^m-1)/(q-1)$. 
Let $\chi_N'$ be a multiplicative character of order $N$ of $\F_{q^m}$ and 
$\chi_N$ be its lift to $\F_{q^{2m}}$.  Furthermore, let $\chi_4$ be  a multiplicative character of order $4$ of $\F_{q^{2m}}$. 
Then, it holds that $G_{q^{2m}}(\chi_4\chi_{N})=G_{q^{2m}}(\chi_4^3\chi_{N})$. 
In particular, it holds that  
\[
G_{q^{2m}}(\chi_4\chi_{N})=\rho_{q^m} \frac{q^m G_{q^m}(\eta {\chi'}_N^{2}) }{G_{q^m}(\eta)}, 
\]
where $\rho_{q^m}=-1$ or $1$ depending on whether $q^m\equiv 3\,(\mod{8})$ or 
 $q^m\equiv 7\,(\mod{8})$, and $\eta$ is the quadratic character of $\F_{q^m}$. 
\end{theorem} 
In \cite{BLMX}, the authors treated only the case where $m=3$ of Theorem~\ref{Gaussmain}. 
The general case can be proved in a similar way.  
\subsection{Gauss periods and Cayley graphs on finite fields}
The {\it Gauss periods} of order $N$ of $\F_{q^m}$ are defined by 
\[
\psi_{\F_{q^m}}(C_i^{(N,q^m)}):=\sum_{x\in C_i^{(N,q^m)}}\psi_{\F_{q^m}}(x), \, \, \, 0\le i\le N-1,
\] 
where $\psi_{\F_{q^m}}$ is the canonical additive character of $\F_{q^m}$. 
By the orthogonality of characters, the Gauss periods of order $N$ of 
$\F_{q^m}$ can be expressed as a linear combination of Gauss sums of $\F_{q^m}$:
\begin{equation}\label{eq:ortho}
\psi_{\F_{q^m}}(C_i^{(N,q^m)})=\frac{1}{N}\sum_{j=0}^{N-1}G_{q^m}(\chi_N^{j})\chi_N^{-j}(\omega^i), \; \, 0\le i\le N-1,
\end{equation}
where $\chi_N$ is a fixed multiplicative character of order $N$ of $\F_{q^m}$ and 
$\omega$ is a fixed primitive element of $\F_{q^m}$.  

\begin{theorem}\label{thm:Yama}{\em (\cite[Theorem~1]{YY})}
Let $\chi$ be a nontrivial multiplicative character of $\F_{q^m}$ and
$\chi'$ be its restriction to $\F_q$.  Take a system $L$ of representatives of $\F_{q^m}^\ast/\F_q^\ast$ such that $\Tr_{q^m/q}$ maps $L$ onto $\{0,1\}\subset\F_q$.  Partition $L$ 
into two parts: 
\[
L_0=\{x\in L:\,\Tr_{q^m/q}(x)=0\} \, \, \mbox{and}\, \, L_1=\{x\in L:\,\Tr_{q^m/q}(x)=1\}, 
\] 
where $\Tr_{q^m/q}$ is the trace function from $\F_{q^m}$ to $\F_q$. 
Then, 
\[
\sum_{x\in L_1}\chi(x)=\left\{
\begin{array}{ll}
G_{q^m}(\chi)/G_q(\chi'), & \mbox{ if   $\chi'$ is nontrivial,}\\
-G_{q^m}(\chi)/q,& \mbox{ otherwise}. 
 \end{array}
\right.
\]
\end{theorem}
Set $S=\{i\,(\mod{\frac{q^m-1}{q-1}}): \omega^i\in L_0\}$, that is, the so-called {\it Singer difference set}. Note that $|S|=(q^{m-1}-1)/(q-1)$.  
If $N\,|\,(q^m-1)/(q-1)$, the restriction  of a multiplicative character $\chi_N$ of order $N$ of $\F_{q^m}$ to $\F_q$ is trivial. In this case, 
by \eqref{eq:ortho} and Theorem~\ref{thm:Yama}, we have 
\begin{align}\label{eq:Gaussperio}
\psi_{\F_{q^m}}(C_i^{(N,q^m)})=&\, -\frac{1}{N}+\frac{1}{N}\sum_{j=1}^{N-1}G_{q^m}(\chi_N^{j})\chi_N^{-j}(\omega^i)\nonumber\\
=&\, -\frac{q^m-1}{N(q-1)}+\frac{q}{N}\sum_{j=0}^{N-1}\sum_{\ell\in S}\chi_N^j(\omega^{\ell-i}). 
\end{align}
Let $\overline{S_N}$ be the reduction  of $S$ modulo $N$. We may identify $\overline{S_N}$ with a group ring element in $\Z[\Z_N]$:  
\[
\overline{S_N}=c_0 [0]+c_1[1]+\cdots +c_{N-1} [N-1]\in \Z[\Z_N], 
\]
where $c_i\in \Z$, $i=0,1,\ldots,N-1$. 
Define 
\begin{equation}\label{eq:defFN}
F_N=\{c_i:0\le i\le N-1\}
\end{equation}
as an ordinary set, i.e., the set of multiplicities of $x\in\Z_{N}$ in the multiset $\overline{S_N}$. Furthermore, 
let 
\begin{equation}\label{eq:defIIII}
I_\beta=\{i\in \Z_N:c_i =\beta\}, \, \beta \in F_N. 
\end{equation}
Then,  
\begin{equation}\label{eq:sNequ}
\overline{S_N}=\sum_{\beta\in F_N}\beta I_\beta \in \Z[\Z_N].
\end{equation} 
Continuing from \eqref{eq:Gaussperio}, we have 
\begin{equation}\label{eq:fiform}
\psi_{\F_{q^m}}(C_i^{(N,q^m)})= -\frac{q^m-1}{N(q-1)}+q\beta,  
\end{equation}
where $\beta\in F_N$ is determined by $i\in I_\beta$. Thus, the values of Gauss periods of order $N$ of $\F_{q^m}$ are computable from $\overline{S_N}$.  

\begin{remark}\label{rem:todo}
In this paper, we will take a union of cyclotomic classes of order  
$4N$ of $\F_{q^{2m}}$ as a connection set of our strongly regular Cayley graph, where $N\,|\,(q^m-1)/(q-1)$. Therefore, we need to compute a sum of Gauss periods of order $4N$ of $\F_{q^{2m}}$. More precisely, if the connection set has the form  $D=\bigcup_{i\in I}C_{i}^{(4N,q^{2m})}$ for some $I\subseteq \Z_{4N}$, we need to show  that 
$\sum_{i\in I}\psi_{\F_{q^{2m}}}(C_{i+j}^{(4N,q^{2m})})$, $j=0,1,\ldots,4N-1$,
take exactly two values. 
\end{remark}
\section{Three-valued Gauss periods}\label{section:infinite}
\subsection{Brief review of three-valued Gauss periods}\label{section:fun}
Schmidt and White~\cite{SW} studied when Gauss periods $\psi_{\F_{q^m}}(C_i^{(N,q^m)})$, $i=0,1,\ldots,N-1$, with $N\,|\,(q^m-1)/(q-1)$ take exactly two values in relation to a classification problem of projective two-weight irreducible cyclic code. They found two families and eleven sporadic examples of two-valued Gauss periods.   In \cite{BLMX,MX2,MX}, we  constructed strongly regular Cayley graphs based on their results. 

As a natural continuation of the study in \cite{SW}, Feng, Xiang and the author~\cite{FMX2} considered Gauss periods which take exactly 
three rational values, and use them to construct circulant weighing matrices 
and association schemes. In particular, the authors treated three-valued Gauss periods forming an arithmetic progression.  

Let $q^m$ be a prime power and $N>2$ be a positive integer dividing $(q^m-1)/(q-1)$.  
Assume that the Gauss periods $\psi_{\F_{q^m}}(C_i^{(N,q^m)})$, $0\le  i \le N-1$, take exactly three  values $\alpha_1,\alpha_2,\alpha_3$, which
form an arithmetic progression, say, $\alpha_1-\alpha_2=\alpha_2-\alpha_3=t>0$. 
By \eqref{eq:fiform} and the definition \eqref{eq:defFN} of $F_N$,  
we have 
\begin{equation}\label{eq:abab}
F_N=\left\{(\beta_i:=)\frac{\alpha_i}{q}+\frac{q^m-1}{qN(q-1)}:i=1,2,3\right\}. 
\end{equation}
It is clear that, $\beta_1,\beta_2$ and $\beta_3$ also form an arithmetic progression. Write $I_i:=I_{\beta_i}$ for simplicity, where $I_{\beta_i}$ is defined in \eqref{eq:defIIII}. 
Then, the sizes of $I_i$'s are determined as follows. 
\begin{lemma}\label{lem:sizeofI}{\em (\cite[Lemma~2.5]{FMX2})} 
With notation as above, it holds that 
\begin{align*}
|I_1|=&\,\frac{N(\alpha_2^2-\alpha_2t+k)+2\alpha_2-k-t+1}{2t^2}, \\
|I_2|=&\,\frac{N(t^2-\alpha_2^2-k)-1-2\alpha_2+k}{t^2}, \\
|I_3|=&\,\frac{N(\alpha_2^2+\alpha_2t+k)+2\alpha_2-k+t+1}{2t^2}, 
\end{align*}
where $k:=(q^m-1)/N$. 
\end{lemma}

We now have a simple representation of Gauss sums. Let $\omega$ be a primitive element of $\F_{q^m}$ and $\chi$ be a nontrivial multiplicative character of exponent  $N$ of $\F_{q^m}$. Then, by Theorem~\ref{thm:Yama} and \eqref{eq:sNequ}, 
\begin{equation}\label{eq:qads}
G_{q^m}(\chi)=q\sum_{i\in \overline{S_N}}\chi(\omega^i)
=q\sum_{j=1,2,3}\beta_j\sum_{i\in I_j}\chi(\omega^i). 
\end{equation}
Furthermore, by \eqref{eq:abab}, the right-hand side of \eqref{eq:qads} is computed as 
\begin{equation}\label{eq:Iexp}
\sum_{j=1,2,3}\alpha_j\sum_{i\in I_j}\chi(\omega^i)
=t\left(2\sum_{i\in I_1}\chi(\omega^i)+\sum_{i\in I_2}\chi(\omega^i)\right). 
\end{equation}

The following two infinite families of three-valued Gauss 
periods of order $N$ of $\F_{q^m}$, which form an arithmetic progression, were known~\cite{FMX2}: 
\begin{align}
& \mbox{ $m=6$ and $N=(q^3-1)/(q-1)$, \hspace{4.3cm}} \label{eq:known1}\\
& \mbox{ $m=3$ and $N=(q^3-1)/(3(q-1))$ with $q\equiv 1\,(\mod{3})$. \hspace{4.3cm}} \label{eq:known2}
\end{align}
On the other hand, many other examples not belonging to the two families above were known, see \cite[Example~4.4, Table~1]{FMX2}. 
\subsection{Two families of three-valued Gauss periods}\label{sec:newfa}
In this subsection, we find two families of three-valued Gauss periods  forming an arithmetic progression. 
Let $q$ be a prime power and $\omega$ be a primitive 
element of $\F_{q^3}$. Let $M$ be a positive integer dividing  $(q^3-1)/(q-1)$,  and set $N=\frac{q^{3}-1}{M(q-1)}$.  
\begin{lemma}\label{lem:three1}
Assume that $\omega^{j_1N},\omega^{j_2N},\omega^{j_3N}$ are linearly independent over $\F_q$ for all distinct $j_1,j_2,j_3\in \{0,1,\ldots,M-1\}$. 
Then, the Gauss periods $\psi_{\F_{q^3}}(C_i^{(N,q^3)})$, $i=0,1,\ldots,N-1$, take exactly three 
values $-M+2q,-M+q,-M$. 
\end{lemma}
\proof 
Let $\omega^a\in \F_{q^3}^\ast$.  It is impossible that $\Tr_{q^{3}/q}(\omega^a\omega^{jN})=0$ for three or more $j\in \{0,1,\ldots,M-1\}$. In fact,  
if
\[
\Tr_{q^{3}/q}(\omega^a\omega^{j_1N})=\Tr_{q^{3}/q}(\omega^a\omega^{j_2N})=\Tr_{q^{3}/q}(\omega^a\omega^{j_3N})=0
\] for some distinct $j_1,j_2,j_3\in \{0,1,\ldots,M-1\}$, since $\omega^{j_1N},\omega^{j_2N},\omega^{j_3N}$ are linearly independent over $\F_q$,   $\Tr_{q^{3}/q}(\omega^ax)=0$ for all $x\in \F_{q^3}$, which is impossible. Hence, 
\begin{align*}
\psi_{\F_{q^3}}(\omega^a C_0^{(N,q^3)})=&\sum_{j=0}^{M-1}
\psi_{\F_{q}}(\Tr_{q^{3}/q}(\omega^a \omega^{jN})\F_{q}^\ast)\\
=&\left\{
\begin{array}{ll}
-M+2q,&  \mbox{if $\Tr_{q^{3}/q}(\omega^a \omega^{jN})=0$ for exactly two $j\in \{0,1,\ldots,M-1\}$, }\\
-M+q,&  \mbox{if $\Tr_{q^{3}/q}(\omega^a \omega^{jN})=0$ for exactly one  $j\in \{0,1,\ldots,M-1\}$, }\\
-M,&  \mbox{if $\Tr_{q^{3}/q}(\omega^a \omega^{jN})\not=0$ for any $j\in \{0,1,\ldots,M-1\}$.}
 \end{array}
\right.
\end{align*}
This completes the proof of the lemma. 
\qed
\vspace{0.3cm}

Put $\alpha_1=-M+2q$, $\alpha_2=-M+q$ and $\alpha_3=-M$. Then, by \eqref{eq:abab}, we have
$\beta_1=2$, $\beta_2=1$ and $\beta_3=0$.  
Define 
\begin{equation}\label{def:II}
I_j:=\{i\,(\mod{N}): 0\le i\le N-1, \psi_{\F_{q^3}}(C_i^{(N,q^3)})=\alpha_j\}, \, \, j=1,2,3.
\end{equation}
By Lemma~\ref{lem:sizeofI}, we have 
\[
|I_1|=\frac{M-1}{2}, |I_2|=q-M+2,|I_3|=\frac{q^2+q+1}{M}-q+\frac{M-3}{2}. 
\]
Furthermore, by \eqref{eq:sNequ},  the 
reduction of the Singer difference set $S=\{i\,(\mod{q^2+q+1}):\Tr_{q^3/q}(\omega^i)=0\}$ modulo $N$ is 
\begin{equation}\label{eq:quotients}
\overline{S_N}=I_1\cup I_1\cup I_2\subseteq \Z_{N}.
\end{equation}
\begin{corollary}
{\em (\cite[Section~4.3]{FMX2})}\label{cor:three1}
Let $q\equiv 1\,(\mod{3})$ be a prime power, and set $N=(q^2+q+1)/3$. Then, 
the Gauss periods $\psi_{\F_{q^3}}(C_i^{(N,q^3)})$, $i=0,1,\ldots,N-1$, take exactly three 
values $-3+2q,-3+q,-3$. 
\end{corollary}
\proof 
The minimal polynomial of $\omega^N$ is of degree $3$, and hence $1,\omega^N,\omega^{2N}$ are 
linearly independent over $\F_{q}$. Then, by Lemma~\ref{lem:three1}, 
the conclusion of the corollary follows.  
\qed

\vspace{0.3cm}
The result above recovers the family \eqref{eq:known2} of known three-valued Gauss periods. 
\begin{corollary}\label{cor:three2}
Let $q$ be an odd prime power such that $q\equiv 2$ or $4\,(\mod{7})$, and set $N=(q^2+q+1)/7$. Then, 
the Gauss periods $\psi_{\F_{q^3}}(C_i^{(N,q^3)})$, $i=0,1,\ldots,N-1$, take exactly three 
values $-7+2q,-7+q,-7$. 
\end{corollary}
\proof 
The minimal polynomial of each $\omega^{jN}$, $j=1,2,\ldots,6$, is of degree $3$,  
and hence $1,\omega^{jN},\omega^{2jN}$ are 
linearly independent over $\F_{q}$. This also implies that 
$1,\omega^{jN},\omega^{j'N}$ are 
linearly independent over $\F_{q}$, where $j'\equiv 2j\,(\mod{7})$, since 
$\omega^{2jN}\F_{q}^\ast=\omega^{j'N}\F_{q}^\ast$. 

We next show that $1,\omega^{jN},\omega^{3jN}$ are 
linearly independent over $\F_{q}$ for every $j=1,2,\ldots,6$.  If this is shown, similarly to the above, $1,\omega^{jN},\omega^{j'N}$ are 
linearly independent over $\F_{q}$, where $j'\equiv 3j\,(\mod{7})$. 
Assume that $1,\omega^{jN},\omega^{3jN}$ are linearly dependent over $\F_{q}$, i.e., there are $a,b\in \F_q^\ast$ such that $\omega^{3jN}+a\omega^{jN}+b=0$. Then, $f(x)=x^3+ax+b\in \F_q[x]$ is the minimal polynomial of $\omega^{jN}$. Since $\omega^{qjN},\omega^{q^2jN}$ are also roots of $f(x)$, we have $x^3+ax+b=(x-\omega^{jN})(x-\omega^{qjN})(x-\omega^{q^2jN})$. Comparing the coefficients of $x^2$ of both sides, we have $\omega^{jN}+\omega^{qjN}+\omega^{q^2jN}=0$. 
Depending on whether $q\equiv 2$ or $4\,(\mod{7})$, we have 
$1+\omega^{j\frac{q^3-1}{7}}+\omega^{3j\frac{q^3-1}{7}}=0$ or $1+\omega^{2j\frac{q^3-1}{7}}+\omega^{3j\frac{q^3-1}{7}}=0$, 
respectively. In the case where $q\equiv 2\,(\mod{7})$, $g(x)=x^3+x+1$ is the minimal polynomial of $\omega^{j\frac{q^3-1}{7}}$. Since $\omega^{2j\frac{q^3-1}{7}},\omega^{4j\frac{q^3-1}{7}}$ are also roots of $g(x)$, we have $x^3+x+1=(x-\omega^{j\frac{q^3-1}{7}})(x-\omega^{2j\frac{q^3-1}{7}})(x-\omega^{4j\frac{q^3-1}{7}})$. Comparing the constants of both sides, we have $1=-\omega^{7j\frac{q^3-1}{7}}=-1$, which is a contradiction.  
The proof for the case where $q\equiv 4\,(\mod{7})$ is similar. 

Note that $1,\omega^{jN},\omega^{j'N}$ are linearly independent over $\F_q$ if and only if 
so are $\omega^{tN},\omega^{(j+t)N}$, $\omega^{(j'+t)N}$ for every  
$t\in \{0,1,\ldots,6\}$, where the exponents $j+t$ and $j'+t$ are reduced modulo $7$. 
Then, it is straightforward to check that  $\omega^{j_1N},\omega^{j_2N},\omega^{j_3N}$ are linearly independent over $\F_{q}$ for all distinct $j_1,j_2,j_3\in \{0,1,\ldots,6\}$. 
Then, by Lemma~\ref{lem:three1}, 
the conclusion of the corollary follows.  
\qed

\begin{remark}
We give a comment for the case where $q$ is even in Corollary~\ref{cor:three2} although we need only the case where $q$ is odd  in this paper. If $q$ is even, $\F_{q^3}$ contains the subfield of order $8$. Then, $\omega^{\frac{q^3-1}{7}}\in \F_8$ is a root  of $f(x)=1+x^j+x^3\in \F_2[x]$ for $j=1$ or $2$. This implies that $1,\omega^{N j_1},\omega^{N j_2}$ are linearly dependent over $\F_{q}$, where 
$j_1$ and $j_2$ are defined by $j_1\equiv j(q-1)\,(\mod{7})$ and $j_2\equiv 3(q-1)\,(\mod{7})$, respectively. Therefore, we can not apply Lemma~\ref{lem:three1}. 
\end{remark}
\section{Construction of Cayley graphs based on three-valued Gauss periods}\label{section:construction}
In this section, we give a construction of Cayley graphs on $(\F_{q^{2m}},+)$ based on three-valued Gauss periods of $\F_{q^m}$. With some suitable partition of the set $I_2$, the Cayley graph will become  a strongly regular graph. 

We begin with an illustrative example of our construction. 
\begin{example}{\em 
Let $q=7$ and $M=3$. 
As shown in Corollary~\ref{cor:three1}, the Gauss periods $\psi_{\F_{7^3}}(C_i^{(19,7^3)})$, $i=0,1,\ldots,18$, take exactly three values $\alpha_1=11,\alpha_2=4$ and $\alpha_3=-3$. 
With a suitable primitive element of $\F_{7^3}$, we have
\begin{align*}
I_1=&\{ 0 \},\\
I_2=&\{8, 10, 12, 13, 15, 18 \},\\
I_3=&\{  1, 2, 3, 4, 5, 6, 7, 9, 11, 14, 16, 17\}.
\end{align*}
Consider the partition $S_1=\{8,12,18\}, S_2=\{10,13,15\}$ of $I_2$, and   
set $S_i''\equiv 4^{-1}S_i\,(\mod{19})$, $i=1,2$, that is, $S_1''=\{2, 3, 14\}$ and $S_2''=\{8,12, 18\}$. Define
\begin{align*}
Y&\,=\{19i+4j \,(\mod{76}): (i,j)\in (\{0,3\}\times S_1'')\cup   (\{1,2\}\times S_2'')\}  \nonumber\\
&\, \, \, \, \cup \{19i+4j\,(\mod{76}): i=0,1,2,3,\, 
j\in 4^{-1}I_1\,(\mod{19})\}\\
&\,=\{8, 12,37, 56,65,69, 10,15,34,51,67,70, 0,19,38,57\}. 
\end{align*}
Then, with a suitable primitive element of $\F_{7^6}$, $\Cay(\F_{7^6},\bigcup_{i \in Y} C_i^{(76,7^{6})})$ forms a strongly regular graph. }
\end{example}
We now give our construction in a general setting. Let $q^m\equiv 3\,(\mod{4})$ be a prime power and $\omega$ be a fixed primitive element of $\F_{q^m}$. 
Furthermore, let $N$ be an odd positive integer dividing $(q^m-1)/(q-1)$, and let  $C_i^{(N,q^m)}=\omega^i\langle \omega^N\rangle$, $i=0,1,\ldots,N-1$. In this section, we assume that the Gauss periods $\psi_{\F_{q^m}}(C_i^{(N,q^m)})$, $i=0,1,\ldots,N-1$, take exactly three rational values $\alpha_1,\alpha_2,\alpha_3$ forming an arithmetic progression, that is, $(t:=)\alpha_1-\alpha_2=\alpha_2-\alpha_3>0$. Write 
\[
I_j:=\{i\,(\mod{N})\,|\,\psi_{\F_{q^m}}(C_i^{(N,q^m)})=\alpha_j\}, \, \, \, j=1,2,3. 
\]
Let $S_1, S_2$ be a partition of  $I_2$, and 
let  
$S_i'\equiv 2^{-1}S_i\,(\mod{N})$ and $S_i''\equiv 2^{-1}S_i'\,(\mod{N})$ for 
$i=1,2$. Define 
\begin{equation}\label{eq:defX1}
X:=2S_1'' \cup (2S_2''+N)\,(\mod{2N})
\end{equation}
and 
\begin{align}
Y_X:=&\, \{Ni+4j \, (\mod{4N}): (i,j)\in (\{0,3\}\times S_1'')\cup   (\{1,2\}\times S_2'')\}  \nonumber\\
&\, \, \, \, \cup \{Ni+4j\, (\mod{4N}): i=0,1,2,3,\, 
j\in 4^{-1}I_1\,(\mod{N})\}. 
\label{eq:defI}
\end{align} 
Clearly,  $X\equiv 2^{-1}I_2\,(\mod{N})$ and $Y_X\equiv I_1\cup I_1\cup  I_1\cup I_1\cup  I_2\cup I_2\,(\mod{N})$ as multisets. 

The set $D_X$ given in the following will become the connection set of a strongly regular Cayley graph on $(\F_{q^{2m}},+)$ for some suitable partition $S_1,S_2$ of $I_2$. 
\begin{proposition}\label{mainconstruction}
With  notation as above, define 
\begin{equation}\label{eq:dddd}
D_X:=\bigcup_{i \in Y_X}C_i^{(4N,q^{2m})},
\end{equation}
where $C_i^{(4N,q^{2m})}=\gamma^i \langle \gamma^{4N}\rangle $, $i=0,1,\ldots,4N-1$, and 
$\gamma$  a fixed primitive element of $\F_{q^{2m}}$ such that 
$\gamma^{q^m+1}=\omega$. For $a\in \Z_{4N}$, define $b\equiv 4^{-1}a\,(\mod{N})$ and $c\equiv 2b\,(\mod{2N})$. 
Then, the nontrivial character values of $D_X$ are given by 
\begin{align}
\psi_{\F_{q^{2m}}}(\gamma^a D_X)=&\frac{\rho_{q^m} \delta_a q^m}{2G_{q^m}(\eta)}\left(2\psi_{\F_{q^m}}(\omega^c\bigcup_{\ell\in X}C_\ell^{(2N,q^m)})-\psi_{\F_{q^m}}
(\omega^c \bigcup_{\ell\in 2^{-1}I_2}C_\ell^{(N,q^m)})\right)\nonumber\\
&\, \hspace{3.0cm}+\frac{(q^m-1)(2|I_1|+|I_2|)}{2N}+\left\{
\begin{array}{ll}
-q^m,&  \mbox{if $a\in I_1 \,(\mod{N})$,}\\
-\frac{q^m}{2},&  \mbox{if $a\in I_2 \,(\mod{N})$, }\\
0,&  \mbox{if $a\in I_3 \,(\mod{N})$. }
 \end{array}
\right.\label{eq:DD-}
\end{align}
where $\rho_{q^m}=1$ or $-1$ depending on
whether $q^m\equiv 7\,(\mod{8})$ or $q^m\equiv 3\,(\mod{8})$, and $\delta_a$ is defined by 
\begin{equation*}
\delta_a=\left\{
\begin{array}{ll}
1,&  \mbox{if $a\equiv 0,1\,(\mod{4})$ and $N\equiv 1\,(\mod{4})$ }\\
&  \mbox{ or $a\equiv 0,3\,(\mod{4})$ and $N\equiv 3\,(\mod{4})$,}\\
-1,&  \mbox{if $a\equiv 2,3\,(\mod{4})$ and $N\equiv 1\,(\mod{4})$}\\
&  \mbox{ or $a\equiv 1,2\,(\mod{4})$ and $N\equiv 3\,(\mod{4})$.}
 \end{array}
\right.
\end{equation*}
Furthermore, 
$\eta$ is the quadratic character of $\F_{q^m}$. 
\end{proposition}
We give a proof of Proposition~\ref{mainconstruction} below. 

Let $\chi_{4}$, $\chi_{N}$ and $\chi_{4N}$ be multiplicative characters of order $4$, $N$ and $4N$ of $\F_{q^{2m}}$, respectively. 
By \eqref{eq:ortho}, we have  
\begin{equation}\label{eq:fi1}
\psi_{\F_{q^{2m}}}(\gamma^a D_X)=\frac{1}{4N}\sum_{h=0}^{4N-1}G_{q^{2m}}(\chi_{4N}^h)\sum_{i\in Y_X}\chi_{4N}^{-h}(\gamma^{a+i}). 
\end{equation}
By the assumption that $N$ is odd, $\chi_{4N}^h$ is uniquely expressed as $\chi_{4}^{h_1}\chi_{N}^{h_2}$ for some  $(h_1,h_2)\in 
\Z_{4}\times \Z_N$. Then, the right hand side of \eqref{eq:fi1} is rewritten as 
\begin{align}
&\,\frac{1}{4N}
\sum_{h_1=0,1,2,3}\sum_{h_2=0}^{N-1}
G_{q^{2m}}(\chi_{4}^{h_1}\chi_{N}^{h_2})
\big(\sum_{j\in \{0,3\}}\sum_{s\in S_1}\chi_{4}^{-h_1}(\gamma^{a+Nj})\chi_{N}^{-h_2}(\gamma^{a+s})\nonumber\\
&\hspace{1.4cm} 
+\sum_{j\in \{1,2\}}\sum_{s\in S_2}\chi_{4}^{-h_1}(\gamma^{a+Nj})\chi_{N}^{-h_2}(\gamma^{a+s})+\sum_{j\in 0,1,2,3}\sum_{s\in I_1}\chi_{4}^{-h_1}(\gamma^{a+Nj})\chi_{N}^{-h_2}(\gamma^{a+s})\big). 
\label{eq:init}
\end{align}
We now compute \eqref{eq:init} by dividing 
it into 
the three partial sums: $P_1,P_2$ and $P_3$, where 
$P_1$ is the contribution of the summands with $h_1=0$, $P_2$ is the contribution of the summands with $h_1=2$, and 
$P_3$ is the contribution of the summands with $h_1=1$ or $3$. Then, 
\begin{equation}\label{eq:p1p2p3}
\psi_{\F_{q^{2m}}}(\gamma^a D)=P_1+P_2+P_3. 
\end{equation}
It is clear that $P_2=0$ since \[
\sum_{j\in \{0,3\}}\chi_4^{-2}(\gamma^{a+Nj})=\sum_{j\in \{1,2\}}\chi_4^{-2}(\gamma^{a+Nj})=\sum_{j\in \{0,1,2,3\}}\chi_4^{-2}(\gamma^{a+Nj})=0. 
\]
We consider the partial sum $P_1$. 
\begin{lemma}\label{lem:res1}
It holds that 
\[
P_1=\frac{(q^m-1)(2|I_1|+|I_2|)}{2N}+\left\{
\begin{array}{ll}
-q^m,&  \mbox{if $a\in I_1 \,(\mod{N})$,}\\
-\frac{q^m}{2},&  \mbox{if $a\in I_2 \,(\mod{N})$, }\\
0,&  \mbox{if $a\in I_3 \,(\mod{N})$. }
 \end{array}
\right.
\]
\end{lemma}
\proof 
We compute 
\begin{equation}
 P_1= 
\frac{1}{2N}\sum_{h_2=0}^{N-1}
G_{q^{2m}}(\chi_{N}^{h_2})
\left(\sum_{s\in I_2}\chi_{N}^{-h_2}(\gamma^{a+ s})+2\sum_{s\in I_1}\chi_{N}^{-h_2}(\gamma^{a+ s})\right). \label{eq:quadpart}
\end{equation}
Let $\chi_N'$ be the multiplicative character of order $N $ of $\F_{q^m}$ such that 
$\chi_N$ is the lift of $\chi_N'$. 
By \eqref{eq:qads} and \eqref{eq:Iexp},   
\[
G_{q^m}({\chi'}_N^{-h_2})=t\sum_{s\in I_1\cup I_1 \cup I_2}{\chi'}_{N}^{-h_2}(\omega^{s})=t\sum_{s\in I_1\cup I_1 \cup I_2}\chi_{N}^{-h_2}(\gamma^{s}). 
\]
On the other hand, by Theorem~\ref{thm:lift}, 
\[
G_{q^{2m}}(\chi_{N}^{h_2})=-G_{q^m}({\chi'}_N^{h_2})^2. 
\]
Continuing from  \eqref{eq:quadpart}, we have 
\begin{align*}
P_1+\frac{2|I_1|+|I_2|}{2N}=&\, -\frac{1}{2Nt}\sum_{h_2=1}^{N-1}
G_{q^m}({\chi'}_N^{h_2})^2
G_{q^m}({\chi'}_N^{-h_2})
{\chi'}_{N}^{-h_2}(\omega^{a})\\
=&\,  -\frac{q^m}{2Nt}\sum_{h_2=1}^{N-1}
G_{q^m}({\chi'}_N^{h_2})
{\chi'}_{N}^{-h_2}(\omega^{a})\\
=&-\frac{q^m}{2N}\sum_{h_2=0}^{N-1}
\sum_{s\in  I_1\cup I_1 \cup I_2}{\chi'}_{N}^{-h_2}(\omega^{-s+a})+\frac{q^m(2|I_1|+|I_2|)}{2N}\\
=&\, \frac{q^m(2|I_1|+|I_2|)}{2N}+\left\{
\begin{array}{ll}
-q^m,&  \mbox{if $a\in I_1 \,(\mod{N})$,}\\
-\frac{q^m}{2},&  \mbox{if $a\in I_2 \,(\mod{N})$, }\\
0,&  \mbox{if $a\in I_3 \,(\mod{N})$. }
 \end{array}
\right.
\end{align*}
Then, the conclusion of the lemma follows. \qed 

\vspace{0.3cm}
Next, we compute the partial sum $P_3$. 
\begin{lemma}\label{lem:res2}
Let 
$b\equiv 4^{-1}a\,(\mod{N})$ and $c\equiv 2b\,(\mod{2N})$. 
Then, it holds that 
\begin{equation}\label{eq:p2}
P_3=\frac{\rho_{q^m} \delta_a q^m}{2G_{q^m}(\eta)}
\left(2\psi_{\F_{q^m}}(\omega^c\bigcup_{\ell\in X}C_\ell^{(2N,q^m)})-\psi_{\F_{q^m}}(\omega^c \bigcup_{\ell\in 2^{-1}I_2}C_\ell^{(N,q^m)})\right), 
\end{equation}
where  $\rho_{q^m}$ and $\delta_a$ are defined in Proposition~\ref{mainconstruction}. 
\end{lemma}
\proof 
It is clear that \[
\sum_{j\in \{0,1,2,3\}}\chi_4^{-h_1}(\gamma^{a+Nj})=0,\, \, \, h=1,3.\]
Hence, we have 
\begin{align}
P_3=&\,\frac{1}{4N}\sum_{h_1=1,3}
\sum_{h_2=0}^{N-1}
G_{q^{2m}}(\chi_{4}^{h_1}\chi_{N}^{h_2})
\big(\sum_{j\in \{0,3\}}\sum_{s\in S_1}\chi_{4}^{-h_1}(\gamma^{a+Nj})\chi_{N}^{-h_2}(\gamma^{a+s})\nonumber\\ 
&\hspace{2cm}+\sum_{j\in \{1,2\}}\sum_{s\in S_2}\chi_{4}^{-h_1}(\gamma^{a+Nj})\chi_{N}^{-h_2}(\gamma^{a+s})\big). 
\end{align}
Since $G_{q^{2m}}(\chi_{4}\chi_{N}^{h_2})=G_{q^{2m}}(\chi_{4}^3\chi_{N}^{h_2})$ by Theorem~2.7, we have
\begin{align}\label{eq:fi2}
P_3=&\,\frac{\delta_a}{2N}
\sum_{h_2=1}^{N-1}
G_{q^{2m}}(\chi_{4}\chi_{N}^{h_2})
\big(\sum_{s\in S_1}\chi_{N}^{-h_2}(\gamma^{a+s})-\sum_{s\in S_2}\chi_{N}^{-h_2}(\gamma^{a+s})\big)\\
&\hspace{0.5cm}+\frac{\delta_a}{2N}
G_{q^{2m}}(\chi_{4})(|S_1|-|S_2|). \label{eq:P33}
\end{align}
Since $G_{q^{2m}}(\chi_4)=\rho_{q^m} q^m$ by Theorem~\ref{thm:semiprim}, the term \eqref{eq:P33} is computed as 
$\delta_a \rho_{q^m} q^m (|S_1|-|S_2|)/2N$.

We compute the right hand side of \eqref{eq:fi2}.  
Applying Theorem~\ref{Gaussmain}, we have 
\begin{align}
&
\sum_{h_2=1}^{N-1}
G_{q^{2m}}(\chi_{4}\chi_{N}^{h_2})
\big(\sum_{s\in S_1}\chi_{N}^{-h_2}(\gamma^{a+s})-\sum_{s\in S_2}\chi_{N}^{-h_2}(\gamma^{a+s})\big)\nonumber\\
=&\, \frac{\rho_{q^m} q^m}{G_{q^m}(\eta)}\sum_{h_2=1}^{N-1}
G_{q^m}(\eta{\chi'}_N^{2h_2})
\big(\sum_{s\in S_1}{\chi'}_{N}^{-h_2}(\omega^{a+s})-\sum_{s\in S_2}{\chi'}_{N}^{-h_2}(\omega^{a+s})\big)\nonumber\\
=&\,\frac{\rho_{q^m} q^m}{G_{q^m}(\eta)}\sum_{h_2=1}^{N-1}
G_{q^m}(\eta{\chi'}_N^{2h_2}){\chi'}_{N}^{-2h_2}(\omega^{2b})
\big(\sum_{s\in S_1'}{\chi'}_{N}^{-2h_2}(\omega^{s})-\sum_{s\in S_2'}{\chi'}_{N}^{-2h_2}(\omega^{s})\big)\nonumber\\
=&\,\frac{\rho_{q^m} q^m}{G_{q^m}(\eta)}\sum_{h_2=1}^{N-1}
G_{q^m}(\eta{\chi'}_N^{2h_2})
\sum_{\ell\in X}\eta{\chi'}_{N}^{-2h_2}(\omega^{\ell+c}). \label{eq:cent}
\end{align}
By \eqref{eq:ortho}, we have 
\begin{align*}
&\, \sum_{h_2=0}^{N-1}
G_{q^m}(\eta{\chi'}_N^{2h_2})
\sum_{\ell\in X}\eta{\chi'}_{N}^{-2h_2}(\omega^{\ell+c})\\
=&\, 2N\psi_{\F_{q^m}}(\omega^c\bigcup_{\ell\in X}C_\ell^{(2N,q^m)})-N\psi_{\F_{q^m}}(\omega^c\bigcup_{\ell\in 2^{-1}I_2}C_\ell^{(N,q^m)}). 
\end{align*}
Hence, the value of \eqref{eq:cent} is computed as 
\begin{align*}
&\, \frac{\rho_{q^m} q^m}{G_{q^m}(\eta)}\left(-G_{q^m}(\eta)\sum_{\ell\in X}\eta(\omega^{\ell+c}) \right.\\
&\hspace{2.5cm} \left.+2N\psi_{\F_{q^m}}(\omega^c\bigcup_{\ell\in X}C_\ell^{(2N,q^m)})-N\psi_{\F_{q^m}}(\omega^c\bigcup_{\ell\in 2^{-1}I_2}C_\ell^{(N,q^m)})\right). 
\end{align*}
Finally, the value of  $\sum_{\ell\in X}\eta(\omega^{\ell+c})$ is computed as 
$|S_1|-|S_2|$.
Summing up,  the conclusion of the lemma follows. \qed 

\vspace{0.3cm}
Now, the conclusion of Proposition~\ref{mainconstruction} follows from 
\eqref{eq:p1p2p3} and 
Lemmas~\ref{lem:res1} and \ref{lem:res2}. 
\begin{remark}\label{eq:rem:thm}
If  $X$ defined in \eqref{eq:defX1} satisfies that   
\begin{align}\label{eq:3-chara}
&\, 2\psi_{\F_{q^m}}(\omega^c \bigcup_{\ell\in X}C_{\ell}^{(2N,q^m)})-
\psi_{\F_{q^m}}(\omega^c \bigcup_{\ell\in 2^{-1}I_2}C_{\ell}^{(N,q^m)})\\
=&\, \left\{
\begin{array}{ll}
\pm G_{q^m}(\eta), & \mbox{ if $c\in 2^{-1}I_2\,(\mod{N})$,}\\
0, & \mbox{ otherwise, }
 \end{array}
\right. \nonumber
\end{align}
substituting \eqref{eq:3-chara} into \eqref{eq:DD-}, 
the set $D_X$ takes exactly  two  nontrivial character values $(q^m-1)(2|I_1|+|I_2|)/2N$ and $-q^m+(q^m-1)(2|I_1|+|I_2|)/2N$. This implies that    $\Cay(\F_{q^{2m}},D_X)$ is strongly regular. 
In particular, it is straightforward to check that it has parameters 
$(q^{2m},r(q^m+1),q^m+r^2-3r,r^2-r)$ with $r=(|I_2|+2|I_1|)(q^m-1)/2N$, which is of negative Latin square type. 
\end{remark}

\section{Finding $X$ satisfying the condition of Remark~\ref{eq:rem:thm}}\label{sec:find}

\subsection{A partition of a conic in $\PG(2,q)$}\label{section:conic}
In this subsection, we briefly review a known ``good'' partition of a conic in 
$\PG(2,q)$ found in \cite{DDMR,FMX1}.  

Let $q$ be an odd prime power 
and $\omega$ be a primitive element of $\F_{q^3}$. 
Viewing $\F_{q^3}$ as a 3-dimensional vector space over $\F_q$, we will use $\F_{q^3}$ as the underlying vector space of $\PG(2,q)$. The points of $\PG(2,q)$ are $\langle \omega^i\rangle:=\omega^i \F_{q}^\ast$,
$0\le i\le N-1$, and the lines of $\PG(2,q)$ are 
\begin{equation}\label{eqn_Lu}
L_c:=\{\langle x\rangle:\,\Tr_{q^3/q}(\omega^c x)=0\},
\end{equation}
where $0\le c< q^2+q+1$. 

Define a quadratic form $Q: \F_{q^3}\rightarrow \F_q$ by $Q(x):=
\tr_{q^3/q}(x^2)$. It is straightforward to check that $Q$ is nondegenerate. 
Therefore, $Q$ defines a conic $\cQ$ in $\PG(2,q)$, which contains $q+1$ points. Consequently each line $L$ of $\PG(2,q)$ meets $\cQ$ in $0$, $1$ or $2$ points. 

Consider the following subset of $\Z_{q^2+q+1}$: 
\begin{equation}\label{eqn_IQ}
W_\cQ:=\{i\, (\mod{q^2+q+1}):Q(\omega^i)=0\}=\{d_0,d_1,\ldots, d_{q}\},
\end{equation}
where the elements are numbered in any unspecified order. 
Then, \[
\cQ=\{\langle \omega^{d_i}\rangle:\,0\le i\le q\}.
\] 
Furthermore, $W_\cQ\equiv 2^{-1}S\, (\mod{q^2+q+1})$, where $S$ is the 
Singer difference set, i.e.,   
$S=\{i\,(\mod{q^2+q+1}):\Tr_{q^3/q}(\omega^i)=0\}$. 

Define $D_1:=\bigcup_{i\in W_\cQ}C_{i}^{(q^2+q+1,q^3)}$, where 
$C_{i}^{(q^2+q+1,q^3)}$ is represented by $\langle \omega^i\rangle$. 
Furthermore, define 
\[
W_s:=\{i\,(\mod{q^2+q+1}):\Tr_{q^3/q}(\omega^{2i})\in C_0^{(2,q)}\}
\] 
and 
\[
W_n:=\{i\,(\mod{q^2+q+1}):\Tr_{q^3/q}(\omega^{2i})\in C_1^{(2,q)}\}.
\] 
\begin{lemma} {\rm (\cite[Equation~(3.5)]{BLMX})} \label{lem:conic2}
The set $D_1$ takes exactly three  nontrivial character values, that is,  
\[
\psi_{\F_{q^3}}(\omega^{c} D_1)=
\left\{
\begin{array}{ll}
-1,&  \mbox{if $c\,(\mod{q^2+q+1})\in W_\cQ$,}\\
-1+\epsilon q,&  \mbox{if 
$c\,(\mod{q^2+q+1})  \in W_s$, }\\
-1-\epsilon q,&  \mbox{if 
$c\,(\mod{q^2+q+1}) \in W_n$,}
 \end{array}
\right.
\]
where 
$\epsilon=1$ or $-1$ depending on whether $q\equiv 1\,(\mod{4})$ or $3\,(\mod{4})$. 
\end{lemma}
We will consider a partition of  $D_1$. 
For $d_0\in W_\cQ$, we define
\begin{equation}\label{eqn_defX00}
{\mathcal X}_{\mathcal Q}:=\{\omega^{d_i}\Tr_{q^3/q}(\omega^{d_0+d_i}):\,1\le i\le q\}\cup\{2  \omega^{d_0}\}
\end{equation}
and 
\begin{equation}\label{eqn_defX}
X_{\mathcal Q}:=\{\log_{\omega}(x)\,(\mod{2(q^2+q+1)}):\, x\in {\mathcal X}_{\mathcal Q}\}.  
\end{equation}
It is clear that $X_{\mathcal Q}\equiv W_\cQ\pmod{q^2+q+1}$. 
We give an important  property of $X_{\mathcal Q}$ below. 
\begin{lemma} {\em (\cite[Lemma~3.4]{FMX1})} \label{lem:inva}
If we use any other $d_i$ in place of $d_0$ in the definition of ${\mathcal X}_{\mathcal Q}$, then the resulting set $X_{\mathcal Q}'$ satisfies that 
$X_{\mathcal Q}'\equiv X_{\mathcal Q}$ or  $X_{\mathcal Q}+(q^2+q+1)\,(\mod{2(q^2+q+1)})$. 
\end{lemma}

The set $X_{\mathcal Q}\subseteq \Z_{2(q^2+q+1)}$ can be expressed as 
\begin{equation}\label{eqn_defX}
X_{\mathcal Q}=2E_1\cup (2E_2+(q^2+q+1))\pmod{2(q^2+q+1)}
\end{equation}
for some $E_1,E_2\subseteq \Z_{q^2+q+1}$ with $|E_1|+|E_2|=q+1$. That is, we are partitioning $X_{\mathcal Q}$ into the {\it even} and {\it odd} parts. 
Then,  
$2(E_1\cup E_2)\equiv W_\cQ\,(\mod{q^2+q+1})$ and $4(E_1\cup E_2)\equiv S\,(\mod{q^2+q+1})$, i.e., 
$X_{\mathcal Q}$ induces partitions of the conic $\cQ$ and the Singer difference set $S$, respectively. 
Consider the following partition of $D_1$: 
\[
D_{1,1}:=\bigcup_{i\in X_{\mathcal Q}}C_{i}^{(2(q^2+q+1),q^3)} \mbox{\,  and \, }  D_{1,2}:=\bigcup_{i\in X_{\mathcal Q}+(q^2+q+1)}C_{i}^{(2(q^2+q+1),q^3)}.  
\] 
\begin{theorem} {\rm(\cite[Theorem 3.7, Remark 3.8]{FMX1}, \cite[Theorem~3.4]{BLMX})} \label{thm:main2}
With notation as above, 
the set $D_{1,1}$ takes exactly four nontrivial character values, that is, 
\[
\psi_{\F_{q^3}}(\omega^c D_{1,1})=\left\{
\begin{array}{ll}
\frac{-1+\epsilon\eta(2) G_{q^3}(\eta)}{2}, & \mbox{if $c\,(\mod{q^2+q+1}) \in W_\cQ$ }\\
 & \mbox{ \, \, and 
$c \,(\mod{2(q^2+q+1)}) \in X_{\mathcal Q}$,}\\
\frac{-1-\epsilon\eta(2) G_{q^3}(\eta)}{2}, & \mbox{if $c\,(\mod{q^2+q+1}) \in W_\cQ$ }\\
 & \mbox{ \, \,  and $c\,(\mod{2(q^2+q+1)}) \in X_{\mathcal Q}+(q^2+q+1)$,}\\
\frac{-1+\epsilon q}{2}, & \mbox{if $c\,(\mod{q^2+q+1}) \in W_s$}, \\
\frac{-1-\epsilon q}{2}, & \mbox{if $c\,(\mod{q^2+q+1}) \in W_n$,}
 \end{array}
\right.
\]
where $\eta$ is the quadratic character of $\F_{q^3}$. 
\end{theorem}
\subsection{Quotients of a partition of a conic in $\PG(2,q)$}\label{sec:quo}
Let $q$ be an odd prime power and $M$ be a positive integer dividing $q^2+q+1$, and set $N=\frac{q^2+q+1}{M}$. 
We assume that 
the Gauss periods $\psi_{\F_{q^3}}(C_i^{(N,q^3)})$, $i=0,1,\ldots,N-1$, take exactly three 
values $-M+2q,-M+q,-M$. In Subsection~\ref{sec:newfa}, we found two examples of such $M$, 
that is, 
$M=3,7$.  

Let $X_{\mathcal Q}$ be the subset of $\Z_{2(q^2+q+1)}$ defined in \eqref{eqn_defX}, and let $W_{\mathcal Q}$ be the set defined in  \eqref{eqn_IQ}. 
Then, $X_{\mathcal Q}\equiv W_{{\mathcal Q}}\equiv 2^{-1}S\,(\mod{q^2+q+1})$, where $S$ is the Singer difference set. The reduction of $X_{\mathcal Q}$ modulo $N$ is 
the multiset $\overline{S_N}=2^{-1}(I_1\cup I_1\cup I_2)$ as seen in \eqref{eq:quotients}. We are now interested in whether the reduction $\overline{X_{\mathcal Q}}$ as a multiset  of $X_{\mathcal Q}$ modulo $2N$ is ``purely'' a subset of $\Z_{2N}$. Here, we say that a multiset defined in a group $G$ is {\it purely} a subset of $G$ if each element in $G$ appears in the multiset at most once.  

Define 
\begin{equation}\label{eq:defXi}
X_i:=\{x\,(\mod{2N}):x \in X_{\mathcal Q},x\,(\mod{N})\in 2^{-1}I_i\},  \, \, i=1,2. 
\end{equation}
Then, $X_1\equiv 2^{-1}(I_1 \cup I_1)\,(\mod{N})$ and $X_2\equiv 2^{-1}I_2\,(\mod{N})$ as multisets. Furthermore, $\overline{X_{\mathcal Q}}=X_1\cup X_2$. 
Note that $X_2$ is purely a subset of $\Z_{2N}$, but
$X_1$ may not be purely a subset of $\Z_{2N}$. 
It is clear that  $\overline{X_{\mathcal Q}}$ is purely a subset of $\Z_{2N}$ if and only if so is $X_1$.  
\begin{figure}[t]
\hspace{4cm} 
\caption{$W_{\mathcal Q}, \overline{X_{\mathcal Q}}$ and $2^{-1}\overline{S_{N}}$ are  obtained as the reductions of $X_{\mathcal Q}$ modulo $q^2+q+1,2N$ and $N$, respectively. Furthermore, $2^{-1}(I_1\cup I_1)$ and $2^{-1}I_2$ are obtained as the reductions of $X_1$ and $X_2$, respectively, modulo $N$.}\label{figfig1}
\vspace{-5cm}
\centering
\begin{center}
\includegraphics[scale=0.55,bb=0 0 1347 477]{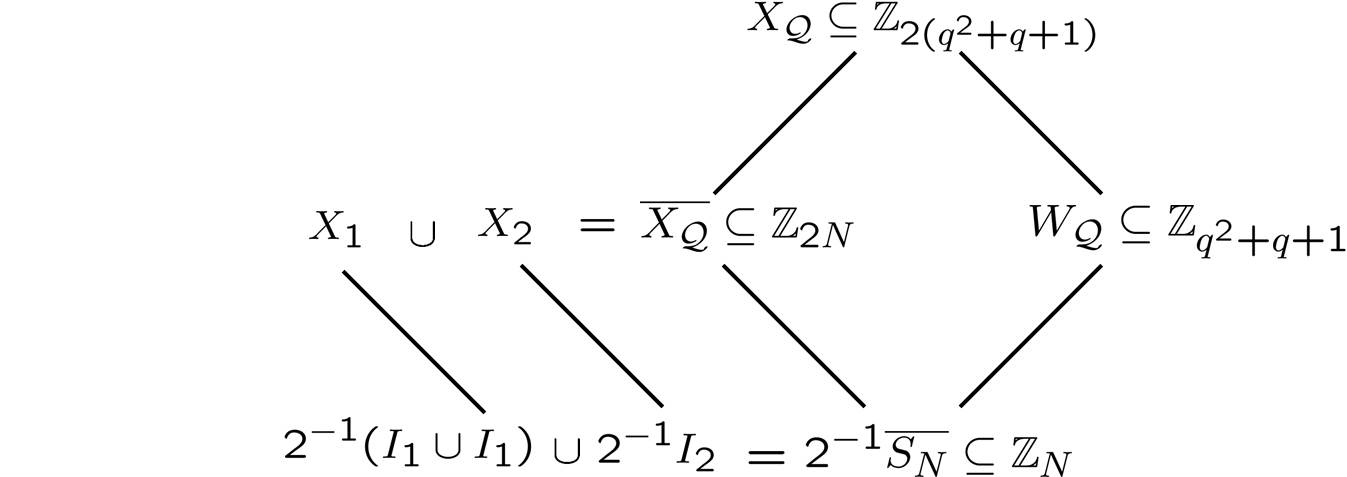}
\end{center}
\end{figure}
The situation is summarized in Figure~\ref{figfig1}. 
\begin{lemma}\label{lem:reduc}
With setting as above, if $X_1$ is purely a subset of $\Z_{2N}$, 
it holds that 
\begin{align*}
&\, 2\psi_{\F_{q^3}}(\omega^c \bigcup_{\ell\in X_2}C_{\ell}^{(2N,q^3)})-
\psi_{\F_{q^3}}(\omega^c \bigcup_{\ell\in 2^{-1}I_2}C_{\ell}^{(N,q^3)})\\
=&\, 2\psi_{\F_{q^3}}(\omega^c \bigcup_{\ell\in \overline{X_{\mathcal Q}}}C_{\ell}^{(2N,q^3)})-
\psi_{\F_{q^3}}(\omega^c \bigcup_{\ell\in 2^{-1}(I_1\cup I_1\cup I_2)}C_{\ell}^{(N,q^3)}). 
\end{align*}
\end{lemma}
\proof
We compute \[
A:=2\psi_{\F_{q^3}}(\omega^c \bigcup_{\ell\in \overline{X_{\mathcal Q}}}C_{\ell}^{(2N,q^3)})-2\psi_{\F_{q^3}}(\omega^c \bigcup_{\ell\in X_2}C_{\ell}^{(2N,q^3)}). 
\]
Since $\overline{X_{\mathcal Q}}=X_1\cup X_2$, we have 
\begin{equation}\label{eq:AA}
A= 2\psi_{\F_{q^3}}(\omega^c \bigcup_{\ell\in X_1}C_{\ell}^{(2N,q^3)}). 
\end{equation}
Since $X_1\equiv 2^{-1}(I_1\cup I_1)\,(\mod{N})$, 
the assumption of the lemma implies that there is a subset $T$ of $X_1$ such that  
\[
X_1=T\cup (T+N)
\]
and 
\[
T\cap (T+N)=\emptyset, 
\]
where $T\equiv 2^{-1}I_1\,(\mod{N})$. Then, continuing from \eqref{eq:AA},  
\begin{align*}
A
=&\, 2\psi_{\F_{q^3}}(\omega^c \bigcup_{\ell\in T}(C_{\ell}^{(2N,q^3)}\cup C_{\ell+N}^{(2N,q^3)}))\\
=&\, 2\psi_{\F_{q^3}}(\omega^c \bigcup_{\ell\in 2^{-1}I_1}C_{\ell}^{(N,q^3)})\\
=&\,\psi_{\F_{q^3}}(\omega^c \bigcup_{\ell\in 2^{-1}(I_1\cup I_1\cup I_2)}C_{\ell}^{(N,q^3)})-\psi_{\F_{q^3}}(\omega^c \bigcup_{\ell\in 2^{-1}I_2}C_{\ell}^{(N,q^3)}).
\end{align*}
This completes the proof of the lemma. 
\qed
\begin{proposition}\label{pro:AB}
Under the assumption of Lemma~\ref{lem:reduc}, it holds that 
\begin{align}
2\psi_{\F_{q^3}}(\omega^c \bigcup_{\ell\in X_2}C_{\ell}^{(2N,q^3)})-
\psi_{\F_{q^3}}(\omega^c \bigcup_{\ell\in 2^{-1}I_2}C_{\ell}^{(N,q^3)})
=\left\{
\begin{array}{ll}
\pm G_{q^3}(\eta), & \mbox{ if $c\in 2^{-1}I_2\,(\mod{N})$,}\\
0, & \mbox{ otherwise. }
 \end{array}
\right. \nonumber
\end{align}
\end{proposition}
\proof
By Lemma~\ref{lem:reduc}, we compute 
\[
B:=2\psi_{\F_{q^3}}(\omega^c \bigcup_{\ell\in \overline{X_{\mathcal Q}}}C_{\ell}^{(2N,q^3)})-
\psi_{\F_{q^3}}(\omega^c \bigcup_{\ell\in 2^{-1}(I_1\cup I_1\cup I_2)}C_{\ell}^{(N,q^3)}). 
\]
It is clear that 
\[
B_1:=2\psi_{\F_{q^3}}(\omega^c \bigcup_{\ell\in \overline{X_{\mathcal Q}}}C_{\ell}^{(2N,q^3)})
= 2\sum_{h=0}^{M-1}\psi_{\F_{q^3}}(\omega^{c+2hN} \bigcup_{\ell\in X_{{\mathcal Q}}}C_{\ell}^{(2(q^2+q+1),q^3)}) 
\]
and 
\[
B_2:=\psi_{\F_{q^3}}(\omega^c \bigcup_{\ell\in 2^{-1}(I_1\cup I_1\cup  I_2)}C_{\ell}^{(N,q^3)})
=\sum_{h=0}^{M-1}\psi_{\F_{q^3}}(\omega^{c+hN} \bigcup_{\ell\in W_{{\mathcal Q}}}C_{\ell}^{(q^2+q+1,q^3)}). 
\]
Then, by Theorem~\ref{thm:main2}, we have 
\begin{align*}
B_1=&\, (-1+\epsilon \eta(2)G_{q^3}(\eta)) \cdot |X_{{\mathcal Q}}\cap (2N\Z_{2(q^2+q+1)}+c)|\\
&\, +(-1-\epsilon \eta(2)G_{q^3}(\eta)) \cdot |(X_{{\mathcal Q}}+(q^2+q+1))\cap (2N\Z_{2(q^2+q+1)}+c)|\\
&\, 
+(-1+\epsilon q)\cdot |W_{s}\cap (N\Z_{q^2+q+1}+c)| 
+(-1-\epsilon q) \cdot |W_{n}\cap (N\Z_{q^2+q+1}+c)|.
\end{align*}
On the other hand, by Lemma~\ref{lem:conic2}, we have 
\begin{align*}
B_2=&\,-|W_{{\mathcal Q}}\cap (N\Z_{q^2+q+1}+c)|\\
&\, +(-1+\epsilon q)\cdot |W_{s}\cap (N\Z_{q^2+q+1}+c)|
+(-1-\epsilon q)\cdot |W_{n}\cap (N\Z_{q^2+q+1}+c)|. 
\end{align*}
Write 
\begin{align*}
x_0:&\,=|X_{{\mathcal Q}}\cap (2N\Z_{2(q^2+q+1)}+c)|,\\ 
x_1:&\,=|(X_{{\mathcal Q}}+(q^2+q+1))\cap (2N\Z_{2(q^2+q+1)}+c)|. 
\end{align*}
Since $x_0+x_1=
|W_{{\mathcal Q}}\cap (N\Z_{q^2+q+1}+c)|$, 
we have 
\begin{equation}\label{eq:BBB}
B=B_1-B_2=\epsilon \eta(2)G_{q^3}(\eta)(x_0-x_1). 
\end{equation}
By the assumption that $X_1$ is purely a subset of $\Z_{2N}$, 
the reduction of $X_{{\mathcal Q}}$ modulo $2N$, that is, $\overline{X_{{\mathcal Q}}}$, is purely a subset of $\Z_{2N}$. Hence, 
$x_0,x_1\in \{0,1\}$. 
In particular, 
\begin{align*}
(x_0,x_1)=(1,1)\, & \Leftrightarrow \, \, c\in 2^{-1}I_1\,(\mod{N}), \\
(x_0,x_1)=(0,1),(1,0)\, & \Leftrightarrow \, \, c\in 2^{-1}I_2\,(\mod{N}),\\
(x_0,x_1)=(0,0)\, & \Leftrightarrow \, \, c\in 2^{-1}I_3\,(\mod{N}). 
\end{align*}
Hence, continuing from \eqref{eq:BBB}, 
we have 
\[
B=\left\{
\begin{array}{ll}
\pm G_{q^3}(\eta), & \mbox{ if $c\in 2^{-1}I_2\,(\mod{N})$,}\\
0, & \mbox{ otherwise. }
 \end{array}
\right. 
\]
This completes the proof of the lemma. 
\qed
\begin{remark}\label{rem:rem1}
Note that $X_2$ is purely a subset of $\Z_{2N}$ and satisfies that $X_2\equiv 2^{-1}I_2\,(\mod{N})$. Then, Proposition~\ref{pro:AB} implies that if  $X_1$  is purely a subset of $\Z_{2N}$, the set  $X_2$ satisfies the condition~\eqref{eq:3-chara} of Remark~\ref{eq:rem:thm} with $m=3$ and $X=X_2$. 
\end{remark}
\subsection{Sufficient conditions for $X_{1}$ to be purely a subset of $\Z_{2N}$}
As shown in Subsection~\ref{sec:newfa}, if  $N=(q^2+q+1)/M$ with $M=3$ or $7$, 
the Gauss periods $\psi_{\F_{q^3}}(C_i^{(N,q^3)})$, $i=0,1,\ldots,N-1$, take exactly three 
values $-M+2q,-M+q,-M$. 
In order to apply 
Remark~\ref{rem:rem1}, we need to investigate when $X_1$ is purely a subset of $\Z_{2N}$. 
The following is our main result of this subsection. 
 \begin{result}\leavevmode\label{result:11}
\begin{itemize}
\item[(i)] 
Let $q$ be a prime power such that $q\equiv 1\,(\mod{3})$, and let $M=3$.
If $q\equiv 7$ or $13\,(\mod{24})$, then $X_1$ is purely  a subset of $\Z_{2N}$. 
\item[(ii)] 
Let $q$ be a prime power such that $q\equiv 2$ or $4\,(\mod{7})$, and let $M=7$.
If $q\equiv 11,37,51$ or $53\,(\mod{56})$, then $X_1$ is purely  a subset of $\Z_{2N}$. 
\end{itemize}
\end{result}

Assume that 
$u\in W_{\mathcal Q}$ and $u\,(\mod{N})\in 2^{-1}I_1$. 
Since $W_{\mathcal Q}\equiv 2^{-1}(I_1\cup I_1\cup I_2)\,(\mod{N})$, there is exactly one $\ell_u\in \{1,2,\ldots,M-1\}$ such that $u+\ell_u N$ is also in $W_{\mathcal Q}$. 
By the definition~\eqref{eqn_IQ} of $W_{\mathcal Q}$, we have 
$\Tr_{q^3/q}(\omega^{2u})=\Tr_{q^3/q}(\omega^{2u+2\ell_u N})=0$. From these equations, 
we have 
\begin{align}
\omega^{2u q^2}&=-(\omega^{2u}+\omega^{2u q}), \label{eq:71}\\
\omega^{2uq}&=-\frac{\omega^{2u-\frac{2\ell_u (q^3-1)}{M}}(1-\omega^{\frac{2(q+1)\ell_u (q^3-1)}{M}})}{1-\omega^{\frac{2q\ell_u (q^3-1)}{M}}}. \label{eq:72}
\end{align}
Define 
\[
g_M(\omega^u)=\Tr_{q^3/q}(\omega^{2u+\ell_u N})
\omega^{\ell_u N} 
\] 
for $u\in W_{\mathcal Q}$ such that $u\,(\mod{N})\in 2^{-1}I_1$. 
By the definition~\eqref{eqn_defX00} of ${\mathcal X}_{\mathcal Q}$ and Lemma~\ref{lem:inva}, we can assume that 
$2\omega^u,g_M(\omega^u) \omega^u\in {\mathcal X}_{\mathcal Q}$ or 
$2\omega^{u+(q^2+q+1)},g_M(\omega^u) \omega^{u+(q^2+q+1)}\in {\mathcal X}_{\mathcal Q}$. 
Thus,  $X_1$ is purely  a subset of $\Z_{2N}$ if and only if  $\eta(2)\not=\eta(g_M(\omega^u))$ for all $u\in W_{\mathcal Q}$ such that   $u\,(\mod{N})\in 2^{-1}I_1$, where $\eta$ is the quadratic character of $\F_{q^3}$. 
\begin{lemma}\label{eq:mm}
With notation as above, 
it holds that 
\[
\eta(g_M(\omega^u))=
\eta(-1)\eta(1-\omega^{\frac{\ell_u (q+1)(q^3-1)}{M}})
\eta(1-\omega^{\frac{2\ell_u q(q^3-1)}{M}}). 
\]
\end{lemma}
\proof 
By the definition of $g_M(\omega^u)$ and the condition~\eqref{eq:71}, we have 
\begin{align}
g_M(\omega^u)&\,=\omega^{2\ell_u N}
(\omega^{2u}+\omega^{2uq+\frac{\ell_u (q^3-1)}{M}}+\omega^{2uq^2+\frac{\ell_u (q+1)(q^3-1)}{M}})\nonumber\\
&\,=\omega^{2\ell_u N}
(\omega^{2u}(1-\omega^{\frac{\ell_u (q+1)(q^3-1)}{M}})+\omega^{2uq+\frac{\ell_u (q^3-1)}{M}}(1-\omega^{\frac{\ell_u q(q^3-1)}{M}})). \label{eq:74}
\end{align}
By substituting \eqref{eq:72} into \eqref{eq:74}, we have 
\begin{equation}\label{eq:eee}
\eqref{eq:74}=-\omega^{2\ell_u N}\omega^{2u-\frac{\ell_u (q^3-1)}{M}}\frac{(1-\omega^{\frac{\ell_u (q^3-1)}{M}})(1-\omega^{\frac{\ell_u (q+1)(q^3-1)}{M}})(1-\omega^{\frac{\ell_u q(q^3-1)}{M}})}{
1-\omega^{\frac{2\ell_u q(q^3-1)}{M}}}. 
\end{equation}
Since $\eta(1-\omega^{\frac{\ell_u (q^3-1)}{M}})=\eta(1-\omega^{\frac{\ell_u q(q^3-1)}{M}})$, 
we have 
\[
\eta(g_M(\omega^u))=\eta(-1)\eta(1-\omega^{\frac{\ell_u (q+1)(q^3-1)}{M}})
\eta(1-\omega^{\frac{2\ell_u q(q^3-1)}{M}}). 
\]
This completes the proof of the lemma. \qed

\begin{proposition}\label{prop:three2}
Let $q\equiv 1\,(\mod{6})$ be a prime power and  
$\eta$ be the quadratic character of $\F_{q^3}$. Then, 
it holds that 
\[
\eta(g_3(\omega^u))=\left\{
\begin{array}{ll}
1, & \mbox{ if $q\equiv 1\,(\mod{12})$,}\\
-1, & \mbox{ if  $q\equiv 7\,(\mod{12})$.}
 \end{array}
\right. 
\]
\end{proposition}
\proof 
By Lemma~\ref{eq:mm}, we have
\begin{align*}
\eta(g_3(\omega^u))=&\,
\eta(-1)\eta(1-\omega^{\frac{2\ell_u (q^3-1)}{3}})
\eta(1-\omega^{\frac{2\ell_u (q^3-1)}{3}})\\
=&\,\eta(-1)=\left\{
\begin{array}{ll}
1, & \mbox{ if $q\equiv 1\,(\mod{12})$,}\\
-1, & \mbox{ if  $q\equiv 7\,(\mod{12})$.}
 \end{array}
\right. 
\end{align*}
This completes the proof of the proposition. 
\qed
\vspace{0.3cm}

By the supplementary low of the quadratic reciprocity, 
$\eta(2)\not=\eta(g_3(\omega^u))$ if and only if 
$q\equiv 7$ or $13\,(\mod{24})$. Thus, we obtain Result~\ref{result:11}~(i). 
\begin{proposition}\label{prop:three2}
Let $q$ be a prime power such that $q\equiv 9$ or $11\,(\mod{14})$ and  
$\eta$ be the quadratic character of $\F_{q^3}$. Then, 
it holds that 
\[
\eta(g_7(\omega^u))=1. 
\]
\end{proposition}
\proof 
By Lemma~\ref{eq:mm}, we have 
\begin{equation}\label{eq:eee2ee}
\eta(g_7(\omega^u))=
\eta(-1)\eta(1-\omega^{\frac{\ell_u (q+1) (q^3-1)}{7}})
\eta(1-\omega^{\frac{2\ell_u q(q^3-1)}{7}}). 
\end{equation}
We consider the right hand side of \eqref{eq:eee2ee} in two cases: 
(i) $q\equiv 9\,(\mod{14})$; and (ii) $q\equiv 11\,(\mod{14})$. 

(i) In the case where $q\equiv 9\,(\mod{14})$, we have 
$1-\omega^{\frac{\ell_u (q+1) (q^3-1)}{7}}=1-\omega^{\frac{3\ell_u (q^3-1)}{7}}$
and 
$1-\omega^{\frac{2\ell_u q (q^3-1)}{7}}=-\omega^{\frac{4\ell_u  (q^3-1)}{7}}(1-\omega^{\frac{3\ell_u  (q^3-1)}{7}})$.  
Then, continuing from \eqref{eq:eee2ee}, we obtain  $\eta(g_7(\omega^u))=1$. 

(ii) In the case where $q\equiv 11\,(\mod{14})$, 
we have 
$1-\omega^{\frac{\ell_u (q+1) (q^3-1)}{7}}=1-\omega^{\frac{5\ell_u (q^3-1)}{7}}$
and 
$1-\omega^{\frac{2\ell_u q (q^3-1)}{7}}=-\omega^{\frac{\ell_u  (q^3-1)}{7}}(1-\omega^{\frac{6\ell_u  (q^3-1)}{7}})$.  Since 
\[
\eta(1-\omega^{\frac{6\ell_u  (q^3-1)}{7}})=\eta(1-\omega^{\frac{6q^2\ell_u  (q^3-1)}{7}})=\eta(1-\omega^{\frac{5\ell_u  (q^3-1)}{7}}), 
\]
continuing from \eqref{eq:eee2ee}, 
we obtain $\eta(g_7(\omega^u))=1$. 
\qed
\vspace{0.3cm}

By the supplementary low of the quadratic reciprocity, 
$\eta(2)=-1$ if and only if 
$q\equiv 3$ or $5\,(\mod{8})$. Thus, we obtain Result~\ref{result:11}~(ii). 

Finally, by applying Proposition~\ref{pro:AB}, Remarks~\ref{eq:rem:thm} and ~\ref{rem:rem1} to Result~\ref{result:11}~(i) and (ii) with the restriction $q^m\equiv 3\,(\mod{4})$, the conclusions of  Theorem~\ref{thm:result}~(i) and (ii) follow, respectively. 
\section{Discussion}\label{section:sporadicexamples}
In this paper, we found two infinite families of strongly regular Cayley graphs on $(\F_{q^{6}},+)$ based on three-valued Gauss periods. 
In this section, we first explain how our study is related to a previous one. We observe that Proposition~\ref{mainconstruction} and Remark~\ref{eq:rem:thm} still hold even in the case where $I_1=\emptyset$;  
in this case, 
the Gauss periods take exactly two values $\alpha_2$ and $\alpha_3$. 
For example, if $N=q^2+q+1$ and $m=3$, we have $I_1=\emptyset,I_2=S,I_3=\Z_N\setminus S$, where $S$ is the Singer difference set.  Furthermore, by Lemma~\ref{lem:conic2} and Theorem~\ref{thm:main2}, $X_{\mathcal Q}$ satisfies that  
\begin{align*}
&\, 2\psi_{\F_{q^3}}(\omega^c \bigcup_{\ell\in X_{\mathcal Q}}C_{\ell}^{(2(q^2+q+1),q^3)})-
\psi_{\F_{q^3}}(\omega^c \bigcup_{\ell\in W_{\mathcal Q}}C_{\ell}^{(q^2+q+1,q^3)})\\
=&\, \left\{
\begin{array}{ll}
\pm G_{q^3}(\eta), & \mbox{ if $c\in W_{\mathcal Q}\,(\mod{q^2+q+1})$,}\\
0, & \mbox{ otherwise, }
 \end{array}
\right. 
\end{align*}
where $W_{\mathcal Q}=2^{-1}I_2$. Hence, by Proposition~\ref{mainconstruction} and Remark~\ref{eq:rem:thm}, we can claim that $\Cay(\F_{q^6},D_{X_{\mathcal Q}})$ forms a strongly regular graph having parameters $(q^{6},r(q^3+1),q^3+r^2-3r,r^2-r)$ with $r=(q^2-1)/2$, where 
$D_{X}$ is defined in \eqref{eq:dddd}. This family of strongly regular graphs was already found in \cite{BLMX}. 
Thus, our construction is a generalization of that given in \cite{BLMX}. Furthermore, for $X_2$ defined in \eqref{eq:defXi}, it is clear that $Y_{X_{2}}$ is obtained as the reduction of $Y_{X_{\mathcal Q}}$ modulo $4N$, where $Y_X$ is defined in \eqref{eq:defI}. Thus, our strongly regular Cayley graphs can be viewed as ``quotients'' of those found in \cite{BLMX}. 

On the other hand, 
we found one sporadic example of a strongly regular Cayley graph based on  Proposition~\ref{mainconstruction}, which is not within the framework above. 
The Gauss periods $\psi_{\F_{7^7}}(C_i^{(29,7^7)})$, $i=0,1,\ldots,28$, 
take exactly three values $\alpha_1=272,\alpha_2=-71,\alpha_3=-414$~\cite[Table1]{FMX2}. With a suitable primitive element of $\F_{7^7}$, we have 
\begin{align*}
I_1=&\{ 8, 10, 12, 15, 18, 26, 27 \},\\
I_2=&\{ 1, 2, 3, 4, 5, 6, 7, 9, 11, 13, 14, 16, 17, 19, 20, 21, 22, 23, 24, 25, 28 \},\\
I_3=&\{ 0\}.
\end{align*}
Define 
$S_1:=I_2$ and $S_2:=\emptyset$, 
which partition the set $I_2$. 
We checked by computer that the set $X\equiv 2S_1'' \cup (2S_2''+29)\,(\mod{58})$ satisfies the condition~\eqref{eq:3-chara} of Remark~\ref{eq:rem:thm}, where 
$S_i''\equiv 4^{-1}S_i\,(\mod{29})$, $i=1,2$. 
Then, by Proposition~\ref{mainconstruction} and Remark~\ref{eq:rem:thm}, $\Cay(\F_{7^{14}},D_X)$ is a strongly regular graph having parameters $(7^{14},r(7^7+1),7^7+r^2-3r,r^2-r)$ with $r=35(7^7-1)/58$. 

We conclude this paper by giving two open problems for future work. 
Although many examples of $(q^m,N)$ which lead to three-valued Gauss periods have been found by computer in \cite[Table1]{FMX2}, only a few of them were theoretically explained in Subsection~\ref{sec:newfa} of this paper. In particular, only the case where $m=3$ and $N=(q^2+q+1)/M$ with $M=7$ was newly characterized. Hence, we give the following problem. 
\begin{problem}
Characterize all pairs $(q,M)$ such that the Gauss periods $\psi_{\F_{q^3}}(C_i^{((q^2+q+1)/M,q^3)})$, $i=0,1,\ldots,(q^2+q+1)/M-1$, take exactly three 
values $-M+2q,-M+q,-M$. 
\end{problem}
Furthermore, in view of  Proposition~\ref{pro:AB} and Remark~\ref{eq:rem:thm}, we give the following problem. 
\begin{problem}
Under the assumption that the Gauss periods $\psi_{\F_{q^3}}(C_i^{((q^2+q+1)/M,q^3)})$, $i=0,1,\ldots,(q^2+q+1)/M-1$, take exactly three 
values $-M+2q,-M+q,-M$, determine when the reduction of $X_{\mathcal Q}$ modulo $2(q^2+q+1)/M$ is purely a subset of $\Z_{2(q^2+q+1)/M}$. 
\end{problem}

\end{document}